\documentclass[11pt]{article}

\usepackage{amsmath}
\usepackage{amssymb}
\usepackage{amsthm}
\usepackage{indentfirst}
\usepackage{cite}

\newtheorem{theorem}{Theorem}[section]
\newtheorem{lemma}{Lemma}[section]
\newtheorem{corollary}{Corollary}[section]
\newtheorem{remark}{Remark}[section]
\newtheorem{example}{Example}

\numberwithin{equation}{section}

\title{\bf BASIC INVARIANTS\\ of Geometric Mappings}
\author{Nenad O. Vesi\'c}
\date{}

\makeatletter
\def\maketag@@@#1{\hbox{\m@th\normalfont\normalsize#1}}
\makeatother

\begin{document}

  \maketitle

\begin{abstract}
    This study is motivated by the researches
    in the field of invariants for geodesic
    and conformal mappings
    presented in \big(T. Y. Thomas, \cite{thomas1}\big) and
    \big(H.
    Weyl, \cite{weyl}\big).
    The Thomas projective parameter and the Weyl projective tensor are
    generalized in this
    article. Generators for vector spaces of invariants of geometric mappings
    are obtained in here.\\[5pt]

    \textbf{Key words:} {invariant of mapping, affine connection,
    Thomas projective parameter, Weyl projective tensor, curvature
    tensor, affine connection}\\[2pt]

    \textbf{Math. Subj. Classification:} {53A55, 53B05, 53C15}
\end{abstract}

\section{Introduction}

An $N$-dimensional manifold $\mathcal M_N$ equipped with an affine
connection $\nabla$ (with torsion) is called the non-symmetric
affine connection space $\mathbb{GA}_N$ \big(see \cite{eisNRG,
mincic1, mincic2, mincic4,
   z4,mica2,mileva1,mica1,mica3}\big).
 As a special case, the manifold $\mathcal M_N$
equipped with a torsion-free affine connection $\overset0\nabla$ is
called the symmetric affine connection space $\mathbb A_N$. More
details about the theory of symmetric affine connection spaces may
be found in \cite{mik5,mik6,mik2,sinjukov}.

   T. Y. Thomas \cite{thomas1} and H.
  Weyl \cite{weyl} started the research about invariants of special diffeomorphisms
  between symmetric affine connection spaces for different applications in physics. Many authors
  have continued the Thomas's and Weyl's works.
   J. Mike\v s
  \cite{mik1, mik2, mik3, mik5, mik6, mik7, mik8, mik10}, I.
  Hinterleitner \cite{mik2, mik5},
  N. S. Sinyukov \cite{sinjukov}, are some of them. Some of invariant
  geometrical object for diffeomorphisms of non-symmetric affine
  connection spaces are obtained in \cite{mica1, mica2, z4, mica3}.

In this paper, as in the previous articles, books and monographs,
the spaces
  $\mathbb{GA}_N$ and $\mathbb{G\overline A}_N$ will be
   the manifold $\mathcal M_N$ equipped with the affine
  connections $\nabla$ and $\overline\nabla=f(\nabla)$.
  A diffeomorphism
  $f:\mathbb{GA}_N\to\mathbb{G\overline A}_N$ which transforms
  the affine connection $\nabla$ of the space
  $\mathbb{GA}_N$ to the affine connection $\overline\nabla$ of
   the space $\mathbb{G\overline A}_N$
  is \emph{the mapping} of the space $\mathbb{GA}_N$.

  In this paper, we will obtain sets of families of invariants
  for diffeomorphisms defined on affine connection spaces with and without torsion. Moreover, we
  will prove that several of these families of invariants are linearly independent.

  \subsection{Affine connection spaces}

  For different applications in physics, for example
  in the Theory of Relativity
  \cite{e1, e2, e3, z2}, affine connection spaces
  with torsion have been studied.

Let $\mathbb{GA}_N$ be a non-symmetric affine connection space. The
affine connection coefficients
   $L^i_{jk}$ of this space are non-symmetric in the indices $j$ and $k$. The
  symmetric and anti-symmetric part of the
   coefficient $L^i_{jk}$ are respectively:

  \begin{eqnarray}
    L^i_{\underline{jk}}=\frac12\big(L^i_{jk}+L^i_{kj}\big)&\mbox{and}&
    L^i_{\underset\vee{jk}}=\frac12\big(L^i_{jk}-L^i_{kj}\big).
  \end{eqnarray}

  The symmetric part $L^i_{\underline{jk}}$
  is the affine connection coefficient
  for a torsion-free affine connection $\overset0\nabla$.
  The manifold $\mathcal M^N$ equipped with the affine connection
  $\overset0\nabla$ \big(whose coefficients are $L^i_{\underline{jk}}$\big) is
  {the associated space $\mathbb A_N$} (of the space
  $\mathbb{GA}_N$).

  The covariant derive of a tensor $a^i_j$ of the type $(1,1)$ with respect to the
  affine connection of the associated space $\mathbb A_N$ is
  (see \cite{sinjukov, mik2, mik5, mik6}):

  \begin{equation}
    a^i_{j|k}=a^i_{j,k}+L^i_{\underline{\alpha k}}a^\alpha_j-
    L^\alpha_{\underline{jk}}a^i_\alpha,
    \label{eq:covderivativesim}
  \end{equation}

  \noindent for the partial
  derivatives $\partial/\partial x^i$ denoted by comma.

  With respect to the affine connection $\overset0\nabla$ and the
  corresponding covariant derivative $|$, one
  Ricci-type identity is obtained. The corresponding
   curvature tensor of the associated space $\mathbb
  A_N$ is (see \cite{sinjukov, mik2, mik5, mik6}):

  \begin{equation}
    R^i_{jmn}=L^i_{\underline{jm},n}-L^i_{\underline{jn},m}+
    L^\alpha_{\underline{jm}}L^i_{\underline{\alpha n}}-
    L^\alpha_{\underline{jn}}L^i_{\underline{\alpha m}}.
    \label{eq:R}
  \end{equation}

  Based on the definitions and results from L. P. Eisenhart \cite{eisNRG, eis01, eis02}, A.
  Einstein \cite{e1,e2,e3},\linebreak M. Prvanovi\'c \cite{mileva1}, S. M.
  Min\v ci\'c defined four kinds of covariant derivaties with
  respect to the affine connection $\nabla$ of the space $\mathbb{GA}_N$
  \cite{mincic1, mincic2, mincic4}:

  \begin{equation}
    \begin{array}{ll}
      a^i_{j\underset1|k}=a^i_{j,k}+L^i_{\alpha
      k}a^\alpha_j-L^\alpha_{jk}a^i_\alpha,&
      a^i_{j\underset2|k}=a^i_{j,k}+L^i_{k\alpha}a^\alpha_j
      -L^\alpha_{kj}a^i_\alpha,\\
      a^i_{j\underset3|k}=a^i_{j,k}+L^i_{\alpha
      k}a^\alpha_j-L^\alpha_{kj}a^i_\alpha,&
      a^i_{j\underset4|k}=a^i_{j,k}+L^i_{k\alpha}a^\alpha_j
      -L^\alpha_{jk}a^i_\alpha.
    \end{array}
    \label{eq:covderivativensim}
  \end{equation}

  With respect to these generalizations of the covariant derivative (\ref{eq:covderivativesim}),\linebreak S. M.
  Min\v ci\'c got four curvature tensors, eight derived
  curvature tensors and fifteen curvature pseudotensors of the space
  $\mathbb{GA}_N$ \cite{mincic1, mincic2, mincic4}.
  Curvature tensors and derived curvature
  tensors of the space $\mathbb{GA}_N$ are elements of the
  family
  \cite{z4}

  \begin{equation}
    K^i_{jmn}=R^i_{jmn}+uL^i_{\underset\vee{jm}|n}+u'L^i_{\underset\vee{jn}|m}+
    vL^\alpha_{\underset\vee{jm}}L^i_{\underset\vee{\alpha n}}+
    v'L^\alpha_{\underset\vee{jn}}L^i_{\underset\vee{\alpha m}}+
    wL^\alpha_{\underset\vee{mn}}L^i_{\underset\vee{\alpha j}},
    \label{eq:K}
  \end{equation}

  \noindent for the curvature tensor $R^i_{jmn}$ of the associated
  space $\mathbb A_N$ and real coefficients $u,u',v,v',w$. Five of
  twelve
  curvature tensors from the family (\ref{eq:K}) are linearly independent
  \cite{mincic2},
  while the rest can be expressed in terms of these five tensors and the curvature tensor
   $R^i_{jmn}$
  of the associated space $\mathbb{A}_N$.

  \bigskip

  Special kind of non-symmetric
  affine connection spaces are
   $N$-dimensional differentiable manifolds equipped with the non-symmetric metric
  tensor $g_{ij}$  of the type $(0,2)$. The symmetric and
  anti-symmetric part of the metric $g_{ij}$ are

  \begin{eqnarray}
    g_{\underline{ij}}=\frac12\big(g_{ij}+g_{ji}\big)&\mbox{and}&
    g_{\underset\vee{ij}}=\frac12\big(g_{ij}-g_{ji}\big).
    \label{eq:gsimantisim}
  \end{eqnarray}

  These spaces are \emph{the
  generalized Riemannian spaces $\mathbb{GR}_N$} (see \cite{eis01, eis02}).
  The affine connection coefficients of the space $\mathbb{GR}_N$
  are the generalized Christoffel symbols of the second kind:

  \begin{equation}
    \Gamma^i_{jk}=\frac12g^{\underline{i\alpha}}\big(g_{j\alpha,k}-
    g_{jk,\alpha}+g_{\alpha k,j}\big),
    \label{eq:GRNChristoffel2nd}
  \end{equation}

  \noindent for
  $\big[g^{\underline{ij}}\big]=\big[g_{\underline{ij}}\big]^{-1}$.
  After symmetrizing the symbols $\Gamma^i_{jk}$ in the indices $j$ and $k$,
  we get that they reduce to the corresponding
  Christoffel symbols  $\Gamma^i_{\underline{jk}}$ obtained from the
  symmetric metric $g_{\underline{ij}}$ of the associated space $\mathbb R_N$.

\section{About invariants}


  Many
  invariants for mappings of torsion-free spaces
  are obtained. Some of them are the Thomas projective parameter,
  the Weyl projective tensor,
  the Weyl conformal curvature, and many others.
  These invariants may be found in the next monographs, books and
  papers:
  Mike\v s \cite{mik1, mik2, mik3, mik5, mik6, mik8, mik7, mik10},
  Sinyukov \cite{sinjukov}, Hinterleitner \cite{mik5}, Berezovski
  \cite{mik7, mik10, mik8, mik5}, etc.

  In this paper, the author's main purpose is to search inceptive
  invariants of different
  mappings defined on affine connection spaces. We will start the
  generalization of the basic invariants in here.
  The next aim of this article is to discover how
  many of the basic invariants are linearly independent.

  Let $f:\mathbb{GA}_N\to\mathbb{G\overline A}_N$ be a
  mapping between non-symmetric affine connection spaces $\mathbb{GA}_N$ and
  $\mathbb{G\overline A}_N$.

  The deformation tensor $P^i_{jk}=
  \overline{L}^i_{jk}-L^i_{jk}$ of this mapping is

  \begin{equation}
    P^i_{jk}=
    \overline\omega{}^i_{jk}-\omega^i_{jk}+\overline\tau{}^i_{jk}-
    \tau^i_{jk},
    \label{eq:P}
  \end{equation}

  \noindent for geometrical objects
  $\omega^i_{jk},\overline\omega{}^i_{jk},\tau^i_{jk},\overline\tau{}^i_{jk}$
  of the type $(1,2)$ such that $\omega^i_{jk}=\omega^i_{kj},$
  $\overline\omega{}^i_{jk}=\overline\omega{}^i_{kj},$
  $\tau^i_{jk}=-\tau^i_{kj},$
  $\overline\tau{}^i_{jk}=-\overline\tau{}^i_{kj}$. If the mapping
  $f$ is equitorsion \cite{mica2, mica3,z4},
 i.e. if $\overline
  L^i_{\underset\vee{jk}}=L^i_{\underset\vee{jk}}$, the equation
  (\ref{eq:P}) reduces to

  \begin{equation}
    P^i_{jk}=\overline\omega{}^i_{jk}-
    \omega^i_{jk}.
    \label{eq:P'}
  \end{equation}

  After symmetrizing the equations (\ref{eq:P}, \ref{eq:P'}) in
  the indices $j$ and $k$, we get

  \begin{equation}
    P^i_{\underline{jk}}=\overline
    \omega{}^i_{jk}-\omega{}^i_{jk}.
    \label{eq:Psim}
  \end{equation}

The deformation tensor $\overline P^i_{{jk}}$ of the inverse map
$f^{-1}:\mathbb{G\overline
  A}_N\to\mathbb{GA}_N$  is

  \begin{equation*}
  \overline P^i_{{jk}}=
  L^i_{jk}-\overline L^i_{jk}=-P^i_{jk}.
  \end{equation*}

   So, the
  following equalities hold:

  \begin{equation}
    P^i_{\underline{jk}}=\overline
    L^i_{\underline{jk}}-L^i_{\underline{jk}}=\overline
    \omega{}^i_{jk}-\omega^i_{jk}=
    \Big(-\frac12\overline P^i_{\underline{jk}}\Big)-
    \Big(-\frac12 P^i_{\underline{jk}}\Big).
    \label{eq:P'sim}
  \end{equation}

   Hence, the next equation is satisfied

  \begin{equation}
    P^i_{\underline{jk}}=
    {\overline\omega}{}^i_{(1).jk}-
    \omega{}^i_{(1).jk}=
    {\overline\omega}{}^i_{(2).jk}-
    \omega{}^i_{(2).jk}=
    {\overline\omega}{}^i_{(3).jk}-
    \omega{}^i_{(3).jk},
    \label{eq:P''sim}
  \end{equation}

  \noindent for the geometrical objects

  \begin{eqnarray}
    \omega{}^i_{(1).jk}=L^i_{\underline{jk}},&
    \omega{}^i_{(2).jk}=\omega^i_{jk},&
    \omega{}^i_{(3).jk}=-\frac12P^i_{\underline{jk}},
    \label{eq:omega}
  \end{eqnarray}

  \noindent and the corresponding
  ${\overline\omega}{}^i_{(1).jk},
  {\overline\omega}{}^i_{(2).jk},
  {\overline\omega}{}^i_{(3).jk}$.

  Here and after, the symbol $(p)$ means
  that the equal-index summation convention does not apply to
  the index $p$.

  After anti-symmetrizing the equation (\ref{eq:P}) in the indices $j$
  and $k$, we get

  \begin{equation}
    P^i_{\underset\vee{jk}}=\xi^i_{jk}=
    \overline L{}^i_{\underset\vee{jk}}-L^i_{\underset\vee{jk}}=
    \overline\tau{}^i_{jk}-\tau^i_{jk}.
    \label{eq:Pantisim}
  \end{equation}

  \subsection{Invariants in symmetric affine connection space}

  With respect to the equation (\ref{eq:P''sim}), we get:

  \begin{eqnarray*}
    {\widetilde{\overline{\mathcal T}}}{}^i_{(1).jk}=
    {\widetilde{\mathcal T}}{}^i_{(1).jk},&
    {\widetilde{\overline{\mathcal T}}}{}^i_{(2).jk}=
    {\widetilde{\mathcal T}}{}^i_{(2).jk},&
    {\widetilde{\overline{\mathcal T}}}{}^i_{(3).jk}=
    {\widetilde{\mathcal T}}{}^i_{(3).jk},
  \end{eqnarray*}

  \noindent for the geometrical objects

  \begin{eqnarray}
    {\widetilde{\mathcal
    T}}{}^i_{(1).jk}=0,&{\widetilde{\mathcal
    T}}{}^i_{(2).jk}=L^i_{\underline{jk}}-\omega^i_{jk},&{\widetilde{\mathcal
    T}}{}^i_{(3).jk}=
    \frac12\big(\overline L^i_{\underline{jk}}+L^i_{\underline{jk}}\big),
    \label{eq:Tbasic}
  \end{eqnarray}

  \noindent and the corresponding
  ${\widetilde{\overline{\mathcal T}}}{}^i_{(1).jk},
  {\widetilde{\overline{\mathcal T}}}{}^i_{(2).jk},
  {\widetilde{\overline{\mathcal T}}}{}^i_{(3).jk}$.

  From the equation

  \begin{equation*}
    \aligned
    {\widetilde{\overline{\mathcal T}}}{}^i_{(p).jm,n}&-
    {\widetilde{\overline{\mathcal T}}}{}^i_{(p).jn,m}+
    {\widetilde{\overline{\mathcal T}}}{}^\alpha_{(p).jm}
    {\widetilde{\overline{\mathcal T}}}{}^i_{(p).\alpha
    n}-{\widetilde{\overline{\mathcal T}}}{}^\alpha_{(p).jn}
    {\widetilde{\overline{\mathcal T}}}{}^i_{(p).\alpha
    m}\\&={\widetilde{{\mathcal T}}}{}^i_{(p).jm,n}-
    {\widetilde{{\mathcal T}}}{}^i_{(p).jn,m}+
    {\widetilde{{\mathcal T}}}{}^\alpha_{(p).jm}
    {\widetilde{{\mathcal T}}}{}^i_{(p).\alpha
    n}-{\widetilde{{\mathcal T}}}{}^\alpha_{(p).jn}
    {\widetilde{{\mathcal T}}}{}^i_{(p).\alpha
    m},
    \endaligned
  \end{equation*}

  \noindent $p=1,2,3$, one obtains that the following equalities
  are satisfied:

  \begin{equation*}
    {\widetilde{\overline{\mathcal W}}}{}^i_{(p).jmn}=
    {\widetilde{\mathcal W}}{}^i_{(p).jmn},
  \end{equation*}

  \noindent for the geometrical objects

  \begin{equation}
    {\widetilde{\mathcal W}}{}^i_{(p).jmn}=
    R^i_{jmn}-\omega{}^i_{(p).jm|n}+
    \omega{}^i_{(p).jn|m}+
    \omega{}^\alpha_{(p).jm}\omega{}^i_{(p).\alpha
    n}-\omega{}^\alpha_{(p).jn}\omega{}^i_{(p).\alpha
    m},
    \label{eq:Wbasic}
  \end{equation}

  \noindent and the corresponding ${\widetilde{\overline{\mathcal
  W}}}{}^i_{(p).jmn}$.

  Therefore, the following lemma holds:

  \begin{lemma}\label{eq:lemaWTsimbasic}
    Let $f:\mathbb{GA}_N\to\mathbb{G\overline A}_N$ be a
    mapping of the associated space $\mathbb A_N$
    characterized by the deformation tensor \emph{(\ref{eq:P''sim})}.

    The geometrical objects ${\widetilde{\mathcal
    T}}{}^i_{(p).jk},p=1,2,3$, given by the equation
    \emph{(\ref{eq:Tbasic})} are invariants of the mapping $f$.

    The geometrical objects ${\widetilde{\mathcal
    W}}{}^i_{(p).jmn},p=1,2,3$, given by the equation
    \emph{(\ref{eq:Wbasic})} are invariants of the mapping $f$.
    \hfill$\Box$
  \end{lemma}

  An invariant ${\widetilde{\mathcal
    T}}{}^i_{(p).jk}$, for $p=1,2,3$, of the mapping
    $f:\mathbb{GA}_N\to\mathbb{G\overline A}_N$ is \emph{the
    \emph(basic\emph) $p$-th class
    associated invariant of the Thomas type}.
    An invariant
    ${\widetilde{\mathcal
    W}}{}^i_{(p).jk}$, for $p=1,2,3$, of the mapping
    $f:\mathbb{GA}_N\to\mathbb{G\overline A}_N$ is \emph{the
    basic $p$-th class associated invariant of the Weyl type}.

    \begin{remark}
      The equalities $${\overline\omega}{}^\alpha_{(3).jm}
      {\overline\omega}{}^i_{(3).\alpha n}=
      \frac14\overline P^\alpha_{\underline{jm}}
      \overline P^i_{\underline{\alpha n}}=
      \frac14P^\alpha_{\underline{jm}}
      P^i_{\underline{\alpha n}}={\omega{}^\alpha_{(3).jm}}
      \omega{}^i_{(3).\alpha n},$$

      \noindent are satisfied. Thus, the invariant
      ${\widetilde{\mathcal W}}{}^i_{(3).jmn}$ reduces to

      \begin{equation}
        {\widetilde{\mathcal W}}{}^i_{(3).jmn}=
        R^i_{jmn}-\omega{}^i_{(3).jm|n}+
        \omega{}^i_{(3).jn|m}.
        \label{eq:W3sim}
      \end{equation}

    \noindent  This invariant is important for researches about
    invariants of mappings characterized by deformation tensors
    $P^i_{\underline{jk}}$ which are not
    expressed in the form (\ref{eq:Psim}). The almost geodesic mappings of
    the first kind are an example of maps such that \big(see
    \cite{mik5, mik7, mik8, mik10}\big).
    \end{remark}

    \begin{corollary}
      Let in the equation \emph{(\ref{eq:omega})} be
      $\omega{}^i_{(2).jk}=\delta^i_j\rho_k+\delta^i_k\rho_j+\sigma_{jk}^i$,
      for a $1$-form $\rho_j$ and a geometrical object
      $\sigma^i_{jk}$ of the type $(1,2)$ symmetric in the indices $j$ and $k$.

      The geometrical objects\emph:

\begin{footnotesize}
      \begin{align}
        &\aligned
        {\widetilde
        T}{}^i_{(2).jk}&=L^i_{\underline{jk}}-\sigma^i_{jk}-
        \frac1{N+1}\Big(\big(L^\alpha_{\underline{j\alpha}}-\sigma^\alpha_{\underline{j\alpha}}\big)
        \delta^i_k+\big(L^\alpha_{\underline{k\alpha}}-\sigma^\alpha_{\underline{k\alpha}}\big)\delta^i_j\Big),
        \endaligned\label{eq:Tderivedsim}\\\displaybreak[0]
        &\aligned
        {\widetilde W}{}^i_{(2).jmn}&=R^i_{jmn}-\sigma^i_{jm|n}+\sigma^i_{jn|m}+
        \sigma^\alpha_{jm}\sigma^i_{\alpha n}-\sigma^\alpha_{jn}
        \sigma^i_{\alpha m}\\&+
        \frac1{N+1}\delta^i_j\big(R_{[mn]}+\sigma^\alpha_{\alpha[m|n]}\big)+\frac
        N{N^2-1}\delta^i_{[m}R_{jn]}+\frac1{N^2-1}\delta^i_{[m}R_{n]j}\\&-\frac1{N^2-1}\delta^i_m\Big(
        \sigma^\alpha_{\alpha[j|n]}+(N+1)\big(
        \sigma^\alpha_{jn|\alpha}-\sigma^\alpha_{j\alpha|n}-
        \sigma^\alpha_{jn}\sigma^\beta_{\alpha\beta}+
        \sigma^\alpha_{j\beta}\sigma^\beta_{n\alpha}
        \big)\Big)\\&+
        \frac1{N^2-1}\delta^i_n\Big(
        \sigma^\alpha_{\alpha[j|m]}+(N+1)\big(\sigma^\alpha_{jm|\alpha}-\sigma^\alpha_{j\alpha|m}-
        \sigma^\alpha_{jm}\sigma^\beta_{\alpha\beta}+\sigma^\alpha_{j\beta}
        \sigma^\beta_{\alpha
        m}\big)\Big)
        \endaligned\label{eq:Wderivedsim}
      \end{align}
\end{footnotesize}

      \noindent are the
      invariants of the mapping $f$ of the Thomas and the Weyl type
      respectively.
    \end{corollary}

    \begin{proof}
    The geometrical objects:

      \begin{equation}
        {\widetilde{\mathcal T}}{}^i_{(2).jk}=
        L^i_{\underline{jk}}-\delta^i_k\rho_j-\delta^i_j\rho_k-\sigma_{jk}^i
        \label{eq:Tsigma2basic}
      \end{equation}

      \noindent and

      \begin{equation}
        \aligned
        {\widetilde{\mathcal W}}{}^i_{(2).jmn}&=R^i_{jmn}-
        \delta^i_j\rho_{[m|n]}-\delta^i_m\big(\rho_{j|n}
        +\rho_j\rho_n+\sigma^\alpha_{jn}\rho_\alpha\big)\\&+
        \delta^i_n\big(\rho_{j|m}+\rho_j\psi_m+\sigma^\alpha_{jm}\rho_\alpha\big)\\&
        -\sigma^i_{jm|n}+\sigma^i_{jn|m}+
        \sigma^\alpha_{jm}\sigma^i_{\alpha
        n}-\sigma^\alpha_{jn}\sigma^i_{\alpha m},
        \endaligned
        \label{eq:Wsigma2basic}
      \end{equation}

      \noindent are the second class basic associated invariants of the mapping
      $f$ of the Thomas and the Weyl type.

      If we contract the identity
      ${\widetilde{\overline{\mathcal T}}}{}^i_{(2).jk}-
      {\widetilde{\mathcal T}}{}^i_{(2).jk}=0$ over
      $i$ and $k$, we get

      \begin{equation}
        \aligned
        (N+1)\big(\overline\rho_j-\rho_j\big)&=\overline
        L^\alpha_{\underline{j\alpha}}-\overline\sigma{}^\alpha_{j\alpha}-
        L^\alpha_{\underline{j\alpha}}+\sigma^\alpha_{j\alpha}.
        \endaligned\label{eq:T-Tsimi=k}
      \end{equation}

      After substituting the equation (\ref{eq:T-Tsimi=k}) into the
      equality ${\widetilde{\overline{\mathcal T}}}{}^i_{(2).jk}-
      {\widetilde{\mathcal T}}{}^i_{(2).jk}=0$, one confirms that
      the following equality holds

      \begin{equation*}
        {\widetilde{\overline T}}{}^i_{(2).jk}=
        {\widetilde T}{}^i_{(2).jk},
      \end{equation*}

      \noindent for ${\widetilde T}{}^i_{(2).jk}$ from
      the equation (\ref{eq:Tderivedsim}) and the corresponding
      ${\widetilde{\overline T}}{}^i_{(2).jk}$.

      Let be

      \begin{eqnarray*}
        \rho_{ij}=\rho_{j|n}+\rho_j\psi_n+\sigma^\alpha_{jn}\rho_\alpha&\mbox{and}&
        \overline{\rho}_{ij}=
        \overline\rho_{j\|n}+\overline\rho_j\overline\rho_n+
        \overline\sigma^\alpha_{jn}\overline\rho_\alpha,
      \end{eqnarray*}

      \noindent for the covariant derivative with respect to the
      affine connection of the torsion-free space $\mathbb{\overline
      A}_N$ denoted by $\|$.

      With respect to this substitution, the equality
      $0={\widetilde{\overline{\mathcal W}}}{}^i_{(2).jmn}-
      {\widetilde{\mathcal W}}{}^i_{(2).jmn}$ transforms to

      \begin{footnotesize}
      \begin{equation}
        \aligned
        0&=\overline
        R^i_{jmn}-R^i_{jmn}-\delta^i_j\big(\overline\rho_{[mn]}-\rho_{[mn]}\big)-
        \delta^i_m\big(\overline\rho_{jn}-\rho_{jn}\big)+
        \delta^i_n\big(\overline\rho_{jm}-\rho_{jm}\big)\\&-
        \overline\sigma^i_{jm\|n}+\overline\sigma^i_{jn\|m}+
        \overline\sigma^\alpha_{jm}\overline\sigma^i_{\alpha n}-
        \overline\sigma^\alpha_{jn}\overline\sigma^i_{\alpha m}+
        \sigma^i_{jm|n}-\sigma^i_{jn|m}-\sigma^\alpha_{jm}\sigma^i_{\alpha
        n}+\sigma^\alpha_{jn}\sigma^i_{\alpha m}.
        \endaligned\label{eq:W2-W2sim=0d1}
      \end{equation}
      \end{footnotesize}

      After contracting
      the equation (\ref{eq:W2-W2sim=0d1}) in the indices
      $i$ and $j$ and using the relations $R^\alpha_{\alpha mn}=
      -\big(R_{mn}-R_{nm}\big)\equiv-R_{[mn]}$,
      for the alternation in the indices $m$ and $n$ denoted by the
      square brackets,
      we get

      \begin{equation*}
        \aligned
        (N+1)\big(\overline\rho_{[mn]}-\rho_{[mn]}\big)&=-\overline
        R_{[mn]}+R_{[mn]}-\overline\sigma^\alpha_{\alpha m\|n}+
        \overline\sigma^\alpha_{\alpha
        n\|m}+
        \sigma^\alpha_{\alpha m|n}-\sigma^\alpha_{\alpha n|m},
        \endaligned
      \end{equation*}

      \noindent i.e.

\begin{footnotesize}
      \begin{equation}
        \aligned
        0&=\overline R^i_{jmn}+\frac1{N+1}\delta^i_j\big(\overline
        R_{[mn]}+\overline\sigma^\alpha_{\alpha[m\|n]}\big)+\delta^i_n\big(\overline
        \rho_{jm}-\rho_{jm}\big)-\delta^i_m\big(\overline\rho_{jn}-\rho_{jn}\big)\\&-
        \overline\sigma^i_{jm\|n}+\overline\sigma^i_{jn\|m}+
        \overline\sigma^\alpha_{jm}\overline\sigma^i_{\alpha n}-
        \overline\sigma^\alpha_{jn}\overline\sigma^i_{\alpha m}-
        R^i_{jmn}-\frac1{N+1}\delta^i_j\big(R_{[mn]}+\sigma^\alpha_{\alpha[m|n]}\big)\\&
        +
        \sigma^i_{jm|n}-\sigma^i_{jn|m}-\sigma^\alpha_{jm}\sigma^i_{\alpha
        n}+\sigma^\alpha_{jn}\sigma^i_{\alpha m}.
        \endaligned\label{eq:W2-W2sim=0d2}
      \end{equation}
\end{footnotesize}

      If we contract the last equation over $i$ and $n$,
      we obtain that the following equation is satisfied

      \begin{equation}
        \aligned
        (N-1)\big(\overline\rho_{jm}-\rho_{jm}\big)&=-\overline
        R_{jm}+\frac1{N+1}\big(\overline
        R_{[jm]}+\overline\sigma^\alpha_{\alpha[j\|m]}\big)\\&
        +R_{jm}-
        \frac1{N+1}\big(R_{[jm]}+\sigma^\alpha_{\alpha[j|m]}\big)\\&+
        \overline\sigma^\alpha_{jm\|\alpha}-\overline\sigma^\alpha_{j\alpha\|m}-
        \overline\sigma^\alpha_{jm}\overline\sigma^\beta_{\alpha\beta}+
        \overline\sigma^\alpha_{j\beta}\overline\sigma^\beta_{\alpha
        m}\\&
        -\sigma^\alpha_{jm|\alpha}+\sigma^\alpha_{j\alpha|m}+
        \sigma^\alpha_{jm}\sigma^\beta_{\alpha\beta}-
        \sigma^\alpha_{j\beta}\sigma^\beta_{\alpha m}.
        \endaligned\label{eq:psijm-psijmd}
      \end{equation}

      Based on the equations (\ref{eq:W2-W2sim=0d2},
      \ref{eq:psijm-psijmd}), one gets

      \begin{equation*}
        {\widetilde{\overline W}}{}^i_{(2).jmn}=
        {\widetilde W}{}^i_{(2).jmn},
      \end{equation*}

      \noindent for the geometrical object ${\widetilde W}{}^i_{(2).jmn}$ from
       the equation (\ref{eq:Wderivedsim}) and the corresponding
      ${\widetilde{\overline
      W}}{}^i_{(2).jmn}$.
    \end{proof}

    The invariants (\ref{eq:Tderivedsim}, \ref{eq:Wderivedsim}) are
    \emph{the \emph(second kind\emph) derived associated invariants of the mapping $f$ of
    the Thomas and Weyl type}, respectively.

          \begin{remark}
        If the deformation tensor $P^i_{jk}$ of a studied
        mapping is expressed
        in the form (\ref{eq:P}) all invariants of the Weyl type reduce to
        the corresponding
        invariants of the second class. If the
        deformation tensor $P^i_{jk}$ of a mapping satisfies a differential
        equation, it may be obtained just
        the Weyl type invariants of the
        third class. The invariants of the first and the second
        class of the Thomas type produce the multiplied families of invariants
        for mappings. Hence, it is enough to obtain invariants $\mathcal W{}^i_{(p).jmn}$
         for one $p$. All other invariants $\mathcal{\widetilde W}$ of the Weyl type
          reduce to the obtained one.
      \end{remark}

    \subsection{Invariants in non-symmetric affine connection space}

    We will generalize the invariants (\ref{eq:Tbasic},
    \ref{eq:Wbasic}) in this part of the paper.

    With respect to the equation (\ref{eq:Pantisim}), one gets

    \begin{equation*}
      \hat{\overline{\mathcal T}}{}^i_{jk}=\hat{\mathcal
      T}{}^i_{jk},
    \end{equation*}

    \noindent for the geometrical object

    \begin{equation}
      \hat{\mathcal T}{}^i_{jk}=L^i_{\underset\vee{jk}}-\tau^i_{jk},
      \label{eq:Tantisim}
    \end{equation}

    \noindent and the corresponding $\hat{\overline{\mathcal
    T}}{}^i_{jk}$. Based on the equations (\ref{eq:Tbasic},
    \ref{eq:Tantisim}), one obtains that it is satisfied the
    equalities

    \begin{eqnarray*}
      {\overline{\mathcal T}}{}^i_{(1).jk}=
      {\mathcal T}{}^i_{(1).jk},&
      {\overline{\mathcal T}}{}^i_{(2).jk}=
      {\mathcal T}{}^i_{(2).jk},&
      {\overline{\mathcal T}}{}^i_{(3).jk}=
      {\mathcal T}{}^i_{(3).jk},
    \end{eqnarray*}

    \noindent for the geometrical objects

\begin{footnotesize}
    \begin{eqnarray}
      {\mathcal
      T}{}^i_{(1).jk}=L^i_{\underset\vee{jk}}-\tau^i_{jk},&
           {\mathcal T}{}^i_{(2).jk}=L^i_{jk}-\omega{}^i_{jk}-\tau^i_{jk},&
                 {\mathcal T}{}^i_{(3).jk}=L^i_{jk}+
                 \frac12P^i_{\underline{jk}}-\tau^i_{jk},
      \label{eq:Tgeneral}
    \end{eqnarray}
\end{footnotesize}

    \noindent and the corresponding
    ${{\mathcal{\overline T}}}{}^i_{(1).jk},
    {{\mathcal{\overline T}}}{}^i_{(2).jk},
    {{\mathcal{\overline T}}}{}^i_{(3).jk}$.

        \begin{lemma}
      Let $f:\mathbb{GA}_N\to\mathbb{G\overline A}_N$ be a mapping
      of the non-symmetric affine connection space
      $\mathbb{GA}_N$. The
      geometrical objects $\hat{\mathcal T}{}^i_{jk}$ and $\mathcal
      T{}^i_{(p).jk}$,
      respectively
      given by the equations \emph{(\ref{eq:Tantisim},
      \ref{eq:Tgeneral})}, are invariants of the mapping $f$.\qed
    \end{lemma}

    \begin{corollary}
      The invariants \emph{(\ref{eq:Tbasic}, \ref{eq:Tantisim},
      \ref{eq:Tgeneral})} satisfy the equation

      \begin{equation}
        {\mathcal T}{}_{(p).jk}^i=\widetilde{\mathcal
        T}{}_{(p).jk}^i+
        \hat{\mathcal T}{}^i_{jk},
        \label{eq:T=T+T}
      \end{equation}

      \noindent for $p=1,2,3$.\qed
    \end{corollary}

    The invariant $\mathcal T^i_{(p).jk}$, for $p=1,2,3$, is \emph{the
    $p$-th class general invariant of the Thomas type}. The invariant
    $\hat{\mathcal T}{}^i_{jk}$ is \emph{the anti-symmetric invariant
    of the Thomas type}.

    From the difference ${\hat{\overline{\mathcal T}}}{}^i_{jm\|n}-
    {\hat{\mathcal T}}{}^i_{jm|n}$ and the equality
    $\hat{\overline{\mathcal T}}{}^\alpha_{jm}
    \hat{\overline{\mathcal T}}{}^i_{\alpha n}-
    \hat{\mathcal T}{}^\alpha_{jm}
    \hat{\mathcal T}{}^i_{\alpha n}=0$, we obtain the following
    transformation rules:

    \begin{align}
      &\aligned
      \overline
      L^i_{\underset\vee{jm}\|n}-
      L^i_{\underset\vee{jm}|n}&=
      \overline\tau^i_{jm\|n}+{\overline\omega}{}^i_{(p_1).\alpha
      n}\big(\overline
      L^\alpha_{\underset\vee{jm}}-\overline\tau^\alpha_{jm}\big)\\&-
      {\overline\omega}{}^\alpha_{(p_2).jn}\big(\overline
      L^i_{\underset\vee{\alpha m}}-\overline\tau^i_{\alpha
      m}\big)-
      {\overline\omega}{}^\alpha_{(p_3).mn}\big(\overline
      L^i_{\underset\vee{j\alpha}}-\overline\tau^i_{j\alpha}\big)\\&-
      \tau^i_{jm|n}-{\omega}{}^i_{(p_1).\alpha
      n}\big(L^\alpha_{\underset\vee{jm}}-\tau^\alpha_{jm}\big)\\&+
      {\omega}{}^\alpha_{(p_2).jn}
      \big(L^i_{\underset\vee{\alpha m}}-\tau^i_{\alpha
      m}\big)+
      {\omega}{}^\alpha_{(p_3).mn}
      \big(L^i_{\underset\vee{j\alpha}}-\tau^i_{j\alpha}\big),
      \endaligned\label{eq:L|toL|}\\\displaybreak[0]
      &\aligned
      \overline L^\alpha_{\underset\vee{jm}}\overline
      L^i_{\underset\vee{\alpha n}}-L^\alpha_{\underset\vee{jm}}
      L^i_{\underset\vee{\alpha n}}&=\overline
      L^\alpha_{\underset\vee{jm}}\overline\tau^i_{\alpha
      n}+\overline L^i_{\underset\vee{\alpha
      n}}\overline\tau^\alpha_{jm}-\overline\tau^\alpha_{jm}\overline\tau^i_{\alpha
      n}\\&-L^\alpha_{\underset\vee{jm}}\tau^i_{\alpha n}-
      L^i_{\underset\vee{\alpha
      n}}\tau^\alpha_{jm}+\tau^\alpha_{jm}\tau^i_{\alpha n},
      \endaligned\label{eq:LLtoLL}
    \end{align}

    \noindent for $p_1,p_2,p_3=1,2$.

    \begin{remark}
      The following equalities are satisfied:

      \begin{equation*}
        \aligned
        {\overline\omega}{}^i_{(3).jk}\hat{\overline{\mathcal
        T}}{}^l_{mn}\!-\!{\omega}{}^i_{(3).jk}\hat{{\mathcal
        T}}{}^l_{mn}&=-\frac12\overline
        P^i_{\underline{jk}}\hat{\overline{\mathcal
        T}}{}^l_{mn}\!+\!\frac12P^i_{\underline{jk}}\hat{{\mathcal
        T}}{}^l_{mn}={\overline\omega}{}^i_{(1).jk}
        \hat{\overline{\mathcal
        T}}{}^l_{mn}\!-\!\omega{}^i_{(1).jk}
        \hat{{\mathcal
        T}}{}^l_{mn}.
        \endaligned
      \end{equation*}

      \noindent Hence, the geometrical objects
      $\omega{}^i_{(1).jk},\omega{}^i_{(2).jk}$, are enough to express
      all transformation rules
      {(\ref{eq:L|toL|})} with respect to a mapping $f:\mathbb{GA}_N\to\mathbb{G
      \overline A}_N$.
    \end{remark}

      Let us currently express the invariants (\ref{eq:Wbasic}) in
      the form

      \begin{equation*}
        {\widetilde{\mathcal W}}{}^i_{(p).jmn}=R^i_{jmn}-
        {\widetilde{\mathcal D}}{}^i_{(p).jmn},
      \end{equation*}

      \noindent for the corresponding geometrical objects
      ${\widetilde{\mathcal D}}{}^i_{(p).jmn}$. From this
      expression, and the equalities ${\widetilde{\overline{\mathcal W}}}{}^i_{(p).jmn}-
      {\widetilde{\mathcal W}}{}^i_{(p).jmn}=0, p=1,2,3$, we get:

      \begin{equation}
        \overline
        R^i_{jmn}-R^i_{jmn}={\widetilde{\overline{\mathcal
        D}}}{}^i_{(p).jmn}-{\widetilde{\mathcal
        D}}{}^i_{(p).jmn}.
        \label{eq:RtoRD}
      \end{equation}

      With respect to the equation (\ref{eq:K}), we obtain that the
      following equation holds

      \begin{equation}
        \aligned
        \overline K^i_{jmn}-K^i_{jmn}&=\overline
        R^i_{jmn}-R^i_{jmn}+u\big(\overline
        L^i_{\underset\vee{jm}\|n}-L^i_{\underset\vee{jm}|n}\big)\\&+
        u'\big(\overline
        L^i_{\underset\vee{jn}\|m}-L^i_{\underset\vee{jn}|m}\big)+
        v\big(\overline L^\alpha_{\underset\vee{jm}}\overline
        L^i_{\underset\vee{\alpha
        n}}-L^\alpha_{\underset\vee{jm}}L^i_{\underset\vee{\alpha
        n}}\big)\\&+v'\big(\overline
        L^\alpha_{\underset\vee{jn}}\overline
        L^i_{\underset\vee{\alpha
        m}}-L^\alpha_{\underset\vee{jn}}L^i_{\underset\vee{\alpha
        m}}\big)+w\big(\overline
        L^\alpha_{\underset\vee{mn}}\overline
        L^i_{\underset\vee{\alpha
        j}}\!-\!L^\alpha_{\underset\vee{mn}}L^i_{\underset\vee{\alpha
        j}}\big).
        \endaligned\label{eq:KtoK}
      \end{equation}

      After substituting the results (\ref{eq:L|toL|},
      \ref{eq:LLtoLL}, \ref{eq:RtoRD}) into the equation
      (\ref{eq:KtoK}), we get

      \begin{equation*}
        {{\overline{\mathcal
        W}}}{}^i_{(p).(p^1).(p^2).jmn}={{{\mathcal
        W}}}{}^i_{(p).(p^1).(p^2).jmn},
      \end{equation*}

      \noindent for $p=1,2,3,p^1_1,\ldots,p^2_3=1,2$, the
      family geometrical objects
\begin{footnotesize}
      \begin{eqnarray}
        \aligned
        {{{\mathcal
        W}}}{}^i_{(p).(p^1).(p^2).jmn}&=K^i_{jmn}-\omega{}^i_{(p).jm|n}+
        \omega{}^i_{(p).jn|m}+
        \omega{}^\alpha_{(p).jm}\omega{}^i_{(p).\alpha
        n}-\omega{}^\alpha_{(p).jn}\omega{}^i_{(p).\alpha
        m}\\&-
        u\Big(\tau^i_{jm|n}+\omega{}^i_{(p^1_1).\alpha
        n}\big(L^\alpha_{\underset\vee{jm}}-\tau^\alpha_{jm}\big)-
        \omega{}^\alpha_{(p^1_2).jn}\big(L^i_{\underset\vee{\alpha
        m}}-\tau^i_{\alpha m}\big)\Big)\\&-
        u'\Big(\tau^i_{jn|m}+\omega{}^i_{(p^2_1).\alpha
        m}\big(L^\alpha_{\underset\vee{jn}}-\tau^\alpha_{jn}\big)-
        \omega{}^\alpha_{(p^2_2).jm}\big(L^i_{\underset\vee{\alpha
        n}}-\tau^i_{\alpha n}\big)\Big)\\&+
        \big(u\omega{}^\alpha_{(p^1_3).mn}+
        u'\omega{}^\alpha_{(p^2_3).mn}\big)
        \big(L^i_{\underset\vee{j\alpha}}-\tau^i_{j\alpha}\big),
        \endaligned\label{eq:Wpp1p2basic}
      \end{eqnarray}
      \end{footnotesize}

      \noindent the real coefficients $u,u',v,v',w$  and the
      corresponding ${\widetilde{\overline{\mathcal
        W}}}{}^i_{(p).(p^1).(p^2).jmn}$.

      The following theorem holds:

      \begin{theorem}
        Let $f:\mathbb{GA}_N\to\mathbb{G\overline A}_N$ be a mapping of
        the non-symmetric affine connection space
        $\mathbb{GA}_N$.  The
        set ${{{\mathcal
        W}}}{}^i_{(p).(p^1).(p^2).jmn},p=1,2,3,
        p^k=(p^k_1,p^k_2,p^k_3), k=1,2,p^k_r=1,2$,
        of families of geometrical objects given in the equation
        \emph{(\ref{eq:Wpp1p2basic})} is the set of families of
        invariants of the Weyl type for the mapping $f$.\hfill\qed
      \end{theorem}

      \begin{corollary}
        The invariants \emph{(\ref{eq:Wbasic},
        \ref{eq:Wpp1p2basic})} of the mapping $f$ satisfy the equation

\begin{footnotesize}
        \begin{equation}
          \aligned
          {{{\mathcal
          W}}}{}^i_{(p).(p^1).(p^2).jmn}&={\widetilde{\mathcal
          W}}{}^i_{(p).jmn}\\&+
        u\Big(L^i_{\underset\vee{jm}|n}-\tau^i_{jm|n}-\omega{}^i_{(p^1_1).\alpha
        n}\big(L^\alpha_{\underset\vee{jm}}-\tau^\alpha_{jm}\big)+
        \omega{}^\alpha_{(p^1_2).jn}\big(L^i_{\underset\vee{\alpha
        m}}-\tau^i_{\alpha m}\big)\Big)\\&+
        u'\Big(L^i_{\underset\vee{jn}|m}-\tau^i_{jn|m}-\omega{}^i_{(p^2_1).\alpha
        m}\big(L^\alpha_{\underset\vee{jn}}-\tau^\alpha_{jn}\big)-
        \omega{}^\alpha_{(p^2_2).jm}\big(L^i_{\underset\vee{\alpha
        n}}-\tau^i_{\alpha n}\big)\Big)\\&+
        \big(u\omega{}^\alpha_{(p^1_3).mn}+
        u'\omega{}^\alpha_{(p^2_3).mn}\big)
        \big(L^i_{\underset\vee{j\alpha}}-\tau^i_{j\alpha}\big)+
        v\big(L^\alpha_{\underset\vee{jm}}-\tau^\alpha_{jm}\big)
        \big(L^i_{\underset\vee{\alpha n}}-\tau^i_{\alpha n}\big)\\&+
        v'\big(L^\alpha_{\underset\vee{jn}}-\tau^\alpha_{jn}\big)
        \big(L^i_{\underset\vee{\alpha m}}-\tau^i_{\alpha m}\big)+
        w\big(L^\alpha_{\underset\vee{mn}}-\tau^\alpha_{mn}\big)
        \big(L^i_{\underset\vee{\alpha j}}-\tau^i_{\alpha j}\big),
          \endaligned\label{eq:WWcor}
        \end{equation}
\end{footnotesize}

        \noindent for $p=1,2,3, p^1_1,\ldots,p^2_3=1,2$.\hfill\qed
      \end{corollary}

      \begin{corollary}
        Let $f:\mathbb{GA}_N\to\mathbb{G\overline A}_N$ be an
        equitorsion mapping between\linebreak non-symmetric affine connection
        spaces $\mathbb{GA}_N$ and $\mathbb{G\overline A}_N$. The
        invariants \emph{(\ref{eq:Tgeneral})} of this mapping reduce
        to the corresponding invariants \emph{(\ref{eq:Tbasic})}. The invariant
        \emph{(\ref{eq:Tantisim})} coincides with the anti-symmetric
        part $L^i_{\underset\vee{jk}}$ of the affine connection
        coefficient $L^i_{jk}$.

        The set
        \emph{(\ref{eq:Wpp1p2basic})} of families of invariants of
        the mapping $f$ reduces to

\begin{footnotesize}
        \begin{equation}
          \aligned
          {{{\mathcal
          W}}}{}^i_{(p).(p^1).(p^2).jmn}&=K^i_{jmn}-\omega{}^i_{(p).jm|n}+
        \omega{}^i_{(p).jn|m}+
        \omega{}^\alpha_{(p).jm}\omega{}^i_{(p).\alpha
        n}-\omega{}^\alpha_{(p).jn}\omega{}^i_{(p).\alpha
        m}\\&-
        u\Big(\omega{}^i_{(p^1_1).\alpha
        n}L^\alpha_{\underset\vee{jm}}-
        \omega{}^\alpha_{(p^1_2).jn}L^i_{\underset\vee{\alpha
        m}}\Big)-
        u'\Big(\omega{}^i_{(p^2_1).\alpha
        m}L^\alpha_{\underset\vee{jn}}-
        \omega{}^\alpha_{(p^2_2).jm}L^i_{\underset\vee{\alpha
        n}}\Big)\\&
        +
        \big(u\omega{}^\alpha_{(p^1_3).mn}+
        u'\omega{}^\alpha_{(p^2_3).mn}\big)
        L^i_{\underset\vee{j\alpha}},
          \endaligned\label{eq:Wpp1p2basice}
        \end{equation}
\end{footnotesize}

        \noindent for $p=1,2,3,p^1_1,\ldots,p^2_3=1,2$.

        The equation \emph{(\ref{eq:WWcor})} reduces to

        \begin{equation}
          \aligned
          {{{\mathcal
          W}}}{}^i_{(p).(p^1).(p^2).jmn}&={\widetilde{\mathcal
          W}}{}^i_{(p).jmn}+
        \big(u\omega{}^\alpha_{(p^1_3).mn}+
        u'\omega{}^\alpha_{(p^2_3).mn}\big)
        L^i_{\underset\vee{j\alpha}}\\&+
        u\Big(L^i_{\underset\vee{jm}|n}-\omega{}^i_{(p^1_1).\alpha
        n}L^\alpha_{\underset\vee{jm}}+
        \omega{}^\alpha_{(p^1_2).jn}L^i_{\underset\vee{\alpha
        m}}\Big)\\&+
        u'\Big(L^i_{\underset\vee{jn}|m}-\omega{}^i_{(p^2_1).\alpha
        m}L^\alpha_{\underset\vee{jn}}-
        \omega{}^\alpha_{(p^2_2).jm}L^i_{\underset\vee{\alpha
        n}}\Big)\\&+
        vL^\alpha_{\underset\vee{jm}}
        L^i_{\underset\vee{\alpha n}}+
        v'L^\alpha_{\underset\vee{jn}}
        L^i_{\underset\vee{\alpha m}}+
        wL^\alpha_{\underset\vee{mn}}L^i_{\underset\vee{\alpha j}},
          \endaligned\label{eq:WWcore}
        \end{equation}

        \noindent for $p=1,2,3,p^1_1,\ldots,p^2_3=1,2$.\hfill\qed
      \end{corollary}

      The sets of invariants $\mathcal W{}^i_{(p).(p^1).(p^2).jmn}$,
      given by the equations (\ref{eq:Wpp1p2basic},
      \ref{eq:Wpp1p2basice}),
      are \emph{the sets of the $(p.p^1.p^2)$-th class
      invariants of the Weyl type}.

      The families of invariants in the set (\ref{eq:Wpp1p2basic}),
      and in the set (\ref{eq:Wpp1p2basice}) as well, are equivalent
      invariants of the mapping
      $f:\mathbb{GA}_N\to\mathbb{G\overline A}_N$. Because these
      families are different in general, we are aimed to find how
      many of these families are linearly independent.

      The set ${{{\mathcal
          W}}}{}^i_{(p).(p^1).(p^2).jmn}$ of invariants may be expressed as

\begin{scriptsize}
      \begin{eqnarray}
        \aligned
        {{{\mathcal
          W}}}{}^i_{(p).(p^1).(p^2).jmn}&={\widetilde{\mathcal
          W}}{}^i_{(p).jmn}+u\hat{\mathcal T}{}^i_{{jm}|n}+u'
          \hat{\mathcal T}{}^i_{{jn}|m}+
          v\hat{\mathcal T}{}^\alpha_{{jm}}
          \hat{\mathcal T}{}^i_{{\alpha n}}+
          v'\hat{\mathcal T}{}^\alpha_{jn}\hat{\mathcal
          T}{}^i_{\alpha m}+w\hat{\mathcal
          T}{}^\alpha_{mn}\hat{\mathcal T}{}^i_{\alpha j}\\&-
          c^1_puL^i_{\underline{\alpha n}}\hat{\mathcal
          T}{}^\alpha_{jm}-c^2_pu\omega{}^i_{(p).\alpha n}\hat{\mathcal
          T}{}^\alpha_{jm}+c^3_puL^\alpha_{\underline{jn}}\hat{\mathcal
          T}{}^i_{\alpha m}+c^4_pu\omega{}^\alpha_{(p).jn}\hat{\mathcal
          T}{}^i_{\alpha m}\\&-
          c^7_pu'L^i_{\underline{\alpha m}}\hat{\mathcal
          T}{}^\alpha_{jn}-c^8_pu'\omega{}^i_{(p).\alpha m}\hat{\mathcal
          T}{}^\alpha_{jn}+c^9_pu'L^\alpha_{\underline{jm}}\hat{\mathcal
          T}{}^i_{\alpha
          n}+c^{10}_pu'\omega{}^\alpha_{(p).jm}\hat{\mathcal
          T}{}^i_{\alpha
          n}\\&+\big(c^5_pu+c^{11}_pu'\big)L^\alpha_{\underline{mn}}
          \hat{\mathcal T}{}^i_{j\alpha}+
          \big(c^6_pu+c^{12}_pu'\big)\omega{}_{(p).mn}\hat{\mathcal
          T}{}^i_{j\alpha},
        \endaligned\label{eq:Wbasicgen64}
      \end{eqnarray}
\end{scriptsize}

      \noindent for the corresponding coefficients
      $c^k_p\in\{0,1\},k=1,\ldots,12$.

      There are $64$ families of invariants in the set
      (\ref{eq:Wbasicgen64}) characterized by the corresponding
       $11$-tuples

\begin{footnotesize}
       \begin{equation*}
         c_k=\big(1,-c^1_{pk}u,-c^2_{pk}u,c^3_{pk}u,c^4_{pk}u,-c^7_{pk}u',-c^8_{pk}u',c^9_{pk}u',c^{10}_{pk}u',
         c^5_{pk}u+c^{11}_{pk}u',c^6_{pk}u+c^{12}_{pk}u'\big),
       \end{equation*}
\end{footnotesize}

       \noindent $k=1,\ldots,64$, for $c^r_{pk}=c^r_p$ in the $k$-th
       of families of invariants in the set ${\widetilde{{\mathcal
          W}}}{}^i_{(p).(p^1).(p^2).jmn}$.

          The rank of the matrix $
            \left[\begin{array}{c}
              c_1\\
              \vdots\\
              c_{64}
            \end{array}\right]$ of the type $64\times13$
            is $6$.

          Hence, the following theorem holds:

          \begin{theorem}\label{thm:spacesofbasicinv}
            The set \emph{(\ref{eq:Wpp1p2basic})}  of invariants of a mapping
            $f:\mathbb{GA}_N\to\mathbb{G\overline A}_N$, generates the
          $6$-dimensional vector space.\qed
          \end{theorem}

    \section{Applications and examples}

    From the above obtained results, we will search
    invariants of equitorsion geodesic mappings defined on a generalized
    Riemannian space $\mathbb{GR}_N$. Furthermore, we
    will obtain an associated basic invariant of Weyl type of an almost
    geodesic mapping \cite{mik5, mik7, mik8, mik10}
    defined on a Riemannian space $\mathbb
    R_N$. These invariants will be applied in the
    examples after theoretical researches.

    \bigskip

    \noindent\textbf{Equitorsion geodesic mappings.} Let
    $f:\mathbb{GR}_N\to\mathbb{G\overline R}_N$ be an equitorsion
    geodesic mapping. This mapping is characterized by the following
    equation

    \begin{equation}
      P^i_{jk}=\psi_j\delta^i_k+\psi_k\delta^i_j,
      \label{eq:gGRNbasic}
    \end{equation}

    \noindent for the deformation tensor $P^i_{jk}=
    \overline L^i_{jk}-L^i_{jk}$ and a $1$-form $\psi_j$.

    After symmetrizing this equation in the indices $j$ and $k$, we get

    \begin{equation}
      P^i_{\underline{jk}}=
      \psi_j\delta^i_k+\psi_k\delta^i_j.
      \label{eq:gGRNbasicsim}
    \end{equation}

     If we contract the last equation over
    $i$ and $k$, we obtain the following expression of the $1$-form
    $\psi_j$:

    \begin{equation*}
      \psi_j=\frac1{N+1}\big(\overline\Gamma^\alpha_{\underline{j\alpha}}-\Gamma^\alpha_{\underline{j\alpha}}\big).
    \end{equation*}

    After substituting this expression into the equation
    (\ref{eq:gGRNbasicsim}) and recalling the equation
    (\ref{eq:Psim}), we get

    \begin{equation}
      \omega{}^i_{(2).jk}=\frac1{N+1}\delta^i_j\Gamma^\alpha_{\underline{k\alpha}}+
      \frac1{N+1}\delta^i_k\Gamma^\alpha_{\underline{j\alpha}}.
      \label{eq:omega2geod}
    \end{equation}

    If we substitute this expression into the equations
    (\ref{eq:Tbasic}, \ref{eq:Wbasic}), we obtain that the geometrical objects:

    \begin{align}
      &\aligned
      {\widetilde{\mathcal
      T}}{}^i_{(2).jk}=\Gamma^i_{\underline{jk}}-\frac1{N+1}\big(\delta^i_j
      \Gamma^\alpha_{\underline{k\alpha}}+\delta^i_k
      \Gamma^\alpha_{\underline{j\alpha}}\big),
      \endaligned\label{eq:Thomasgeod}\\\displaybreak[0]
      &\aligned
      {\widetilde{\mathcal W}}{}^i_{(2).jmn}&=R^i_{jmn}-
      \frac1{(N+1)^2}\delta^i_m\big((N+1)\Gamma^\alpha_{\underline{j\alpha}|n}+
      \Gamma^\alpha_{\underline{j\alpha}}\Gamma^\alpha_{\underline{n\alpha}}\big)\\&+
      \frac1{(N+1)^2}\delta^i_n\big((N+1)\Gamma^\alpha_{\underline{j\alpha}|m}+
      \Gamma^\alpha_{\underline{j\alpha}}\Gamma^\alpha_{\underline{m\alpha}}\big),
      \endaligned\label{eq:Wsimbasicgeod}
    \end{align}

    \noindent are the basic associated invariants of the Thomas and
    the Weyl
    type of the mapping $f$.

    \begin{example}
      Let $\mathbb{GR}_3$ be a generalized Riemannian space
      equipped with the non-symmetric metric

      \begin{equation}
        g_{ij}=\left[\begin{array}{ccc}
          \big(x^1\big)^2&{x^1}&x^2\\
          -x^1&\big(x^2\big)^2&x^3\\
          -x^2&-x^3&\big(x^3\big)^2
        \end{array}\right].
        \label{eq:G}
      \end{equation}

      The symmetric and anti-symmetric part of this
      metric are

      \begin{eqnarray}
        g_{\underline{ij}}=\left[\begin{array}{ccc}
          \big(x^1\big)^2&0&0\\
          0&\big(x^2\big)^2&0\\
          0&0&\big(x^3\big)^2
        \end{array}\right]&\mbox{and}&
        g_{\underset\vee{ij}}=\left[\begin{array}{ccc}
          0&x^1&x^2\\
          -x^1&0&x^3\\
          -x^2&-x^3&0
        \end{array}\right].
        \label{eq:gh}
      \end{eqnarray}

       The contravariant metric tensor $\big[g^{\underline{ij}}\big]
      =\big[g_{\underline{ij}}\big]^{-1}$
      of the space $\mathbb{GR}_3$ is

      \begin{equation}
        g^{\underline{ij}}=\left[\begin{array}{ccc}
          \big(x^1\big)^{-2}&0&0\\
          0&\big(x^2\big)^{-2}&0\\
          0&0&\big(x^3\big)^{-2}
        \end{array}\right].
        \label{eq:gcontravariant}
      \end{equation}

      To obtain the invariants {(\ref{eq:WWcore})}, we need to
      search the corresponding invariant {(\ref{eq:Wbasic})},
      the geometrical objects $\omega{}^i_{(1).jk}$ and
      $\omega{}^i_{(2).jk}$ and the anti-symmetric part
      $\Gamma^i_{\underset\vee{jk}}$ of the generalized
      Christoffel symbol $\Gamma^i_{jk}$ given by the equation (\ref{eq:GRNChristoffel2nd}).

      The affine connection coefficients of the associated
      space $\mathbb R_3$ are the Chrisfoffel symbols:

      \begin{equation}
        \Gamma^{i}_{\underline{jk}}=\omega{}^i_{(1).jk}=\left\{\begin{array}{cl}
        \big(x^i\big)^{-1},&i=j=k\\
        0,&\mbox{otherwise.}
        \end{array}\right.
        \label{eq:Christoffel1st}
      \end{equation}

      \pagebreak

      Furthermore, based on the equation {(\ref{eq:omega2geod})} we get

      \begin{equation}
        \omega{}^i_{(2).jk}=
        \frac1{N+1}\Big(\delta^i_j\big(x^k\big)^{-1}+
        \delta^i_k\big(x^j\big)^{-1}\Big).
        \label{eq:omegaexm}
      \end{equation}

      For this reason, the curvature tensor of the associated space
      $\mathbb R_3$
      is $R^i_{jmn}=0$. Moreover, the associated invariant of Weyl type
      {(\ref{eq:Wsimbasicgeod})} is

      \begin{equation}
      \aligned
        {\widetilde{\mathcal W}}{}^i_{(2).jmn}&=-
        \frac1{(N+1)^2}\delta^i_m\left((N+1)\Big(\big(x^j\big)^{-1}\Big)_{|n}+
        \big(x^j\big)^{-1}\big(x^n\big)^{-1}\right)\\&+
        \frac1{(N+1)^2}\delta^i_n\left((N+1)\Big(\big(x^j\big)^{-1}\Big)_{|m}+
        \big(x^j\big)^{-1}\big(x^m\big)^{-1}\right).
      \endaligned\label{eq:W2simgeodexample}
      \end{equation}

      Because
      $\Gamma^i_{\underset\vee{jk}}=\dfrac12g^{\underline{i\alpha}}
      \big(g_{\underset\vee{j\alpha},k}-
      g_{\underset\vee{jk},\alpha}+g_{\underset\vee{\alpha k},j}\big)$, we
      get

\begin{footnotesize}
      \begin{eqnarray}
        \Gamma^1_{\underset\vee{23}}=\frac12\big(x^1\big)^{-2}=-\Gamma^1_{\underset\vee{32}},&
        \Gamma^2_{\underset\vee{31}}=\frac12\big(x^2\big)^{-2}=-\Gamma^2_{\underset\vee{13}},&
        \Gamma^3_{\underset\vee{12}}=\frac12\big(x^3\big)^{-2}=-\Gamma^3_{\underset\vee{21}},
        \label{eq:Gammaantisim}
      \end{eqnarray}
\end{footnotesize}

      \noindent and $\Gamma^i_{\underset\vee{jk}}=0$ in all other
      cases.

    After substituting the expression {(\ref{eq:Christoffel1st})}
    into the equation
    {(\ref{eq:WWcore})} and with respect to the equations {(\ref{eq:omegaexm} -
    \ref{eq:Gammaantisim})}, for $p^1_1,\ldots,p^2_3\in\{1,2\}$, we get
    that the set of invariants of
    Weyl type of the mapping $f$  is

\begin{footnotesize}
    \begin{equation}
      \aligned
      {\widetilde{\mathcal
      W}}{}^i_{(2).(p^1).(p^2).jmn}&=-
        \frac1{(N+1)^2}\delta^i_m\left((N+1)\Big(\big(x^j\big)^{-1}\Big)_{|n}+
        \big(x^j\big)^{-1}\big(x^n\big)^{-1}\right)\\&+
        \frac1{(N+1)^2}\delta^i_n\left((N+1)\Big(\big(x^j\big)^{-1}\Big)_{|m}+
        \big(x^j\big)^{-1}\big(x^m\big)^{-1}\right)\\&+
        \big(u\omega{}^\alpha_{(p^1_3).mn}+
        u'\omega{}^\alpha_{(p^2_3).mn}\big)\\&+
        u\Big(\Gamma^i_{\underset\vee{jm}|n}-\omega{}^i_{(p^1_1).\alpha
        n}\Gamma^\alpha_{\underset\vee{jm}}+
        \omega{}^\alpha_{(p^1_2).jn}\Gamma^i_{\underset\vee{\alpha
        m}}\Big)\\&+
        u'\Big(\Gamma^i_{\underset\vee{jn}|m}-\omega{}^i_{(p^2_1).\alpha
        m}\Gamma^\alpha_{\underset\vee{jn}}-
        \omega{}^\alpha_{(p^2_2).jm}\Gamma^i_{\underset\vee{\alpha
        n}}\Big)\\&+
        v\Gamma^\alpha_{\underset\vee{jm}}
        \Gamma^i_{\underset\vee{\alpha n}}+
        v'\Gamma^\alpha_{\underset\vee{jn}}
        \Gamma^i_{\underset\vee{\alpha m}}+
        w\Gamma^\alpha_{\underset\vee{mn}}\Gamma^i_{\underset\vee{\alpha j}}
        \Gamma^i_{\underset\vee{j\alpha}}.
      \endaligned
    \end{equation}
\end{footnotesize}
    \end{example}

    \bigskip

    \noindent\textbf{Almost geodesic mappings.} N. S. Sinyukov
    \cite{sinjukov}, J. Mike\v s \cite{mik1, mik5, mik7, mik8, mik10}
    and many other authors have
    developed the concept of geodesics. We will search
    an associated basic invariant of the mapping $f:\mathbb
    R_N\to\mathbb{\overline R}_N$ characterized by the
    equation

    \begin{equation}
      P^i_{\underline{nm}|j}+P^i_{\underline{jm}|n}+
      P^\alpha_{\underline{jm}}P^i_{\underline{\alpha n}}+
      P^\alpha_{\underline{nm}}P^i_{\underline{\alpha j}}=
      \delta^i_{j}a_{mn}+\delta^i_{n}a_{mj},
      \label{eq:omega3mn}
    \end{equation}

    \noindent for the tensor $a_{ij}$ of the type $(0,2)$ symmetric by
    the indices $i$ and $j$. This mapping is called \emph{the almost
    geodesic mapping} of the type $\tilde\pi{}_1$.

    After replacing the indices $m\leftrightarrow n$, we get

    \begin{equation}
      P^i_{\underline{jn}|m}+
      P^i_{\underline{mn}|j}+
      P^\alpha_{\underline{jn}}P^i_{\underline{\alpha m}}
      +P^\alpha_{\underline{mn}}P^i_{\underline{\alpha j}}=
      \delta^i_{j}a_{nm}+
      \delta^i_{m}a_{nj}.
      \label{eq:omega3nm}
    \end{equation}

    If we subtract the equations
     (\ref{eq:omega3mn}) and (\ref{eq:omega3nm}),
    we will obtain that the following equation is satisfied:

    \begin{equation}
      P^i_{\underline{jm}|n}-P^i_{\underline{jn}|m}=
      -P^\alpha_{\underline{jm}}P^i_{\underline{\alpha n}}+
      P^\alpha_{\underline{jn}}P^i_{\underline{\alpha m}}+
      \delta^i_na_{jm}-\delta^i_ma_{jn}.
      \label{eq:pi1P-P}
    \end{equation}

    Based on the equations (\ref{eq:W3sim}, \ref{eq:pi1P-P})
    and the invariance \linebreak$\overline P^\alpha_{\underline{jm}}
    \overline P^i_{\underline{\alpha n}}=P^\alpha_{\underline{jm}}
    P^i_{\underline{\alpha n}}$,
    we obtain that the geometrical object

    \begin{equation}
      {\widetilde{\mathcal W}}{}^i_{(3).jmn}=
      R^i_{jmn}+\frac12\delta^i_ma_{jn}-\frac12\delta^i_na_{jm},
      \label{eq:pi1basic}
    \end{equation}

    \noindent is the associated basic invariant of the Weyl type of the
    mapping $f$.

    After contracting the equality ${\widetilde{\overline{\mathcal
    W}}}{}^i_{(3).jmn}-{\widetilde{\mathcal W}}{}^i_{(3).jmn}=0$ by
    the indices $i$ and $n$, we obtain that the curvature tensor
    $R^i_{jmn}$ is the derived invariant of this mapping.

    \begin{example}
      Let $\mathbb R_3$ be a Riemannian space equipped with the
      symmetric metric $g_{ij}$ given in the equation
      {(\ref{eq:gh})}. Let also $f:\mathbb
      R_3\to\mathbb{\overline R}_3$ be an almost geodesic mapping of
      the type $\tilde\pi_1$.

      As we obtained in the previous example, the curvature tensor
      of the space $\mathbb R_3$ is $R^i_{jmn}=0$. For this reason,
      the geometrical objects

      \begin{eqnarray}
        {\widetilde{\mathcal
        W}}{}^i_{(3).jmn}=\delta^i_ma_{jn}-\delta^i_na_{jm}&\mbox{and}&
        \widetilde W^i_{(3).jmn}=0,
        \label{eq:pi1basic}
      \end{eqnarray}

      \noindent are the associated basic and the associated derived
       invariant of the Weyl type of the
      mapping $f$.
    \end{example}

\section*{Acknowledgement} 
This paper is financially supported by Serbian Ministry of
Education, Science and Technological Developments, Grant No. 174012.

The author expresses his gratitude to referees for their time
dedicated to revision of this paper.

  \bigskip
  \bigskip
  \bigskip

\noindent\textbf{Author:}\medskip

Nenad O. Vesi\'c

Faculty of Science and Mathematics, Ni\v s, Serbia

Department of Mathematics

Serbian Ministry of Education, Grant No. 174012

\textbf{contact e-mail:} n.o.vesic@outlook.com

\end{document}